\documentclass[12pt]{amsart}

\usepackage{epsf,amssymb,times}

\newtheorem{thm}{Theorem}[section]
\newtheorem{prop}[thm]{Proposition}
\newtheorem{defn}[thm]{Definition}
\newtheorem{lemma}[thm]{Lemma}
\newtheorem{cor}[thm]{Corollary}

\newtheorem{q}[thm]{Question}
\newtheorem*{claim}{Claim}
\newtheorem*{rmk}{Remark}
\newtheorem{HW}{HW}
\newtheorem*{oq}{Open Question}

\newcommand{\C}{\mbox{\bf{C}}}
\newcommand{\R}{\mbox{\bf{R}}}
\newcommand{\Z}{\mbox{\bf{Z}}}
\newcommand{\Q}{\mbox{\bf{Q}}}

\renewcommand{\H}{\mbox{\bf H}}
\newcommand{\bdry}{\partial}
\newcommand{\sa}{\rightsquigarrow}
\newcommand{\s}{\vskip.12in}
\newcommand{\n}{\noindent}

\newcommand{\be}{\begin{enumerate}}
\newcommand{\ee}{\end{enumerate}}

\begin{document}

%%%%%%%%%%%%%%%%%%%%%%%%%%%%%%%%%%%%%%%%%

\title{3-Dimensional Methods in Contact Geometry}

\author{Ko Honda}

\date{December 31, 2003.}

%%%%%%%%%%%%%%%%%%%%%%%%%%%%%%%%%%%%%%%%%
                                         
\maketitle

\begin{tableofcontents}
\end{tableofcontents}

A {\em contact manifold} $(M,\xi)$ is a $(2n+1)$-dimensional manifold $M$
equipped with a smooth maximally nonintegrable hyperplane field
$\xi\subset TM$, i.e., locally $\xi=\ker \alpha$, where $\alpha$ is a
1-form which satisfies $\alpha\wedge (d\alpha)^n\not=0$.  Since
$d\alpha$ is a nondegenerate 2-form when restricted to $\xi$, contact
geometry is customarily viewed as the odd-dimensional sibling of
symplectic geometry.  Although contact 
geometry in dimensions $\geq 5$ is still in an incipient state,
contact structures in dimension $3$ are much better understood,
largely due to the fact that symplectic geometry in two dimensions is
just the study of area.  The goal of this article is to explain some
of the recent developments in 3-dimensional contact geometry, with an
emphasis on methods from 3-dimensional topology.  Basic references
include \cite{Ae,El2,Et1,Ge}.  The article \cite{Kz} is similar in spirit
to ours.

Three-dimensional contact geometry lies at the interface between
3- and 4-manifold geometries, and has been an essential part of the flurry 
in low-dimensional geometry and topology over the last 20 years.  In
dimension 3, it relates to foliation theory and knot theory; in
dimension 4, there are rich interactions with symplectic geometry.  In both
dimensions, there are relations with gauge theories such as Seiberg-Witten
theory and Heegaard Floer homology.

\s\n
{\em Acknowledgements.}  This manuscript grew out of a lecture series
given at the Winter School in Contact Geometry in M\"unchen in February
2003 and a minicourse given at the
{\em Geometry and Foliations 2003} conference, held at Ryokoku University in
Kyoto in September 2003.  I would like to thank Kai Cieliebak
and Dieter Kotschick for the former, and Takashi Tsuboi for the
latter, as well as for his hospitality during my visit to the University
of Tokyo and the Tokyo Institute of Technology during the summer and
fall of 2003.   Much of the actual writing took place during this visit.

\section{Introduction}

From now on we will restrict our attention to contact structures on
$3$-manifolds.  We will implicitly assume that our contact
structures $\xi$ on $M$ satisfy the following:   
\be
\item $\xi$ is oriented, and hence given as the kernel of a global
  1-form $\alpha$.   
\item $\alpha\wedge d\alpha>0$, i.e., the contact structure is {\em
  positive}.   
\ee                                
Such contact structures are often said to be {\em cooriented}.

\begin{HW}
Show that if $\xi$ is a smooth oriented $2$-plane field, then $\xi$
can be written as the kernel of a global 1-form $\alpha$. 
\end{HW}

\subsection{First examples}   $\mbox{}$                                  

\begin{figure}[ht]
\epsfxsize=2.5in 
\centerline{\epsfbox{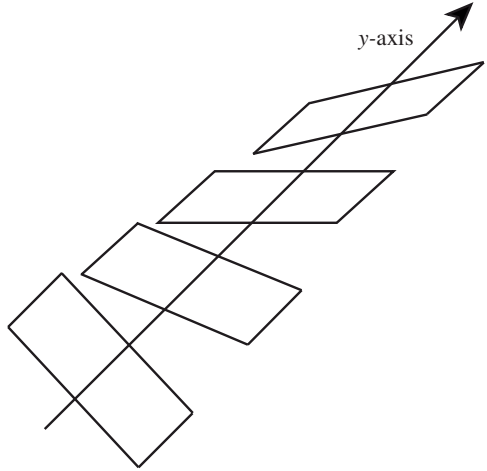}}
\caption{The propeller picture.} 
\label{propeller}	
\end{figure}

\s\n
{\bf Example 1:} $(\R^3,\xi_0)$, where $\R^3$ has coordinates
$(x,y,z)$, and $\xi_0$ is given by $\alpha_0=dz-ydx$.  Then
$\xi_0=\ker \alpha_0=\R\{{\bdry \over \bdry y}, {\bdry \over \bdry x}
+y{\bdry \over \bdry z}\}$.  According to the standard ``propeller
picture'' (see Figure~\ref{propeller}), all the straight lines
parallel to the $y$-axis are 
everywhere tangent to $\xi_0$, and the 2-planes rotate in unison along
these straight lines. 

\s\n
{\bf Example 2:} $(T^3,\xi_n)$.  Here $T^3\simeq \R^3/\Z^3$, with
coordinates $(x,y,z)$, and $n\in \Z^+$.  Then $\xi_n$ is given by
$\alpha_n= \sin(2\pi nz)dx + \cos(2\pi n z)dy$. We have 
$$\xi_n= \R\left\{{\bdry\over \bdry z}, \cos(2\pi n z){\bdry\over
  \bdry x}-\sin(2\pi nz){\bdry \over\bdry y}\right\}.$$  
This time, the circles $x=y=const$ (parallel to the $z$-axis) are
  everywhere tangent to $\xi_n$, and the contact structure makes $n$
  full twists along such circles. 

\begin{HW}
Verify that $(\R^3,\xi_0)$ and $(T^3,\xi_n)$ are indeed contact manifolds.
\end{HW}

The significance of Example 1 is the following:

\begin{thm}[Pfaff] \label{Pfaff}
Every contact 3-manifold $(M,\xi)$ locally looks like $(\R^3,\xi_0)$, i.e., 
for all $p\in M$ there is an open set $U\supset p$ such that
$(U,\xi)\simeq (\R^3,\xi_0)$.  
\end{thm}

Note that an isomorphism in the contact category (usually called a {\em
  contactomorphism}) is a diffeomorphism $\phi: (M_1,\xi_1) \stackrel
  \sim\rightarrow (M_2,\xi_2)$ which maps $\phi_*\xi_1=\xi_2$.  
Pfaff's theorem says that there are no local invariants in contact
  geometry.   

\begin{rmk}
A contactmorphism usually does not preserve the contact 1-form.
\end{rmk}

\begin{HW} 
Prove Pfaff's theorem in dimension 3.  Then generalize it to higher
dimensions. 
\end{HW}

\n
{\bf Example 3:} $(S^3,\xi)$, the standard contact structure on $S^3$.
Consider $B^4=\{|z_1|^2+|z_2|^2\leq 1\} \subset \C^2$.  Then take $S^3=\bdry
B^4$.  The contact structure $\xi$ is defined as follows: for all
$p\in S^3$, $\xi_p$ is the unique complex line $\subset T_pS^3$ (the
unique 2-plane invariant under the complex structure $J$).   

\begin{HW}
Write down a contact 1-form $\alpha$ for $(S^3,\xi)$ and verify that
$\alpha\wedge d\alpha>0$.   
\end{HW}

\subsection{Legendrian knots}
Given a contact manifold $(M,\xi)$, a curve $L\subset M$ is {\em 
Legendrian} if $L$ is everywhere tangent to $\xi$, i.e.,
$\dot L(p)\in \xi_p$ at every point $p\in L$.  In this
section we describe the invariants that can be assigned to a 
Legendrian knot (= embedded closed curve) $L$.  For a more
thorough discussion, see the survey article \cite{Et3}.

\s\n
{\bf Twisting number/Thurston-Bennequin invariant:} Our first invariant is the
{\em relative Thurston-Bennequin invariant} $t(L,\mathcal{F})$,
also known as the {\em twisting number}, where $\mathcal{F}$ is some
fixed framing for $L$.  Although $t(L,\mathcal{F})$ is an invariant of
the {\em unoriented} knot $L$, for convenience pick one orientation of
$L$.  $L$ has a natural framing
called the {\em normal framing}, induced from $\xi$ by taking $v_p\in
\xi_p$ so that $(v_p,\dot L(p))$ form an oriented basis for
$\xi_p$.  We then define $t(L,\mathcal{F})$ to be the integer
difference in the number of twists between the normal framing and
$\mathcal{F}$.  By convention, left twists are negative.  Now, the
framing $\mathcal{F}$ that we choose is often dictated by the
topology.  For example, if $[L]=0\in H_1(M;\Z)$ (which is the
case when $M=S^3$), then there is a compact surface 
$\Sigma\subset M$ with $\bdry \Sigma=L$, i.e., a {\em Seifert
  surface}.  Now $\Sigma$ induces a framing $\mathcal{F}_\Sigma$,
which is the normal framing to the 2-plane field $T\Sigma$ along
$L$, and the {\em Thurston-Bennequin invariant} $tb(L)$ is
given by:  
$$tb(L)=t(L,\mathcal{F}_\Sigma).$$

\begin{HW}
Show that $tb(L)$ does not depend on the choice of Seifert
surface $\Sigma$. 
\end{HW}

In Example 2, if $L=\{x=y=const\}$, then a convenient framing
$\mathcal{F}$ is induced from tori $x=const$ (or equivalently from
$y=const$).  We have $t(L,\mathcal{F})=-n$.

\s\n
{\bf Rotation number:}
Given an oriented Legendrian knot $L$ in $S^3$, we define the
{\em rotation number} $r(L)$ as follows:  Choose a Seifert 
surface $\Sigma$ and trivialize $\xi|_\Sigma$.  Then $r(L)$ is the winding
number of $\dot L$ along $L$ with respect to the
trivialization.

\begin{HW}
Show that $r(L)$ does not depend on the choice of trivialization
or Seifert surface.
\end{HW}

\n
{\bf Front projection:}
We now consider Legendrian knots in the standard contact $(\R^3,\xi_0)$ given
by $dz-ydx=0$.  Consider the {\em front projection} $\pi:\R^3\rightarrow
\R^2$, where $(x,y,z)\mapsto (x,z)$.  Generic Legendrian knots
$L$ (the genericity can be achieved by applying a small contact
isotopy) can be projected to closed curves in $\R^2$ with cusps and ordinary
double points but no vertical tangencies. Conversely, such a closed curve
in $\R^2$ can be lifted to a Legendrian 
knot in $\R^3$ by setting $y$ to be the slope of the curve at $(x,z)$.
(Observe that if $dz-ydx=0$, then $ {dz\over dx}=y$.) The
Thurston-Bennequin invariant and rotation number of a Legendrian knot
$L$ can be computed in the front projection using the following formula:

\begin{eqnarray*}
tb(L) &= & -{1\over 2}(\# \mbox{cusps}) + \# 
\mbox{positive crossings}\\ &&\mbox{ } \mbox{ } -\#\mbox{negative 
crossings}.\\ r(L) &= &{1\over 2} (\#\mbox{downward cusps} 
-\#\mbox{upward cusps})
\end{eqnarray*}

\s
\begin{HW}
Prove the above formulas for $tb$ and $r$ in the front projection.
\end{HW}

\begin{figure}[ht]
\epsfysize=1.7in 
\centerline{\epsfbox{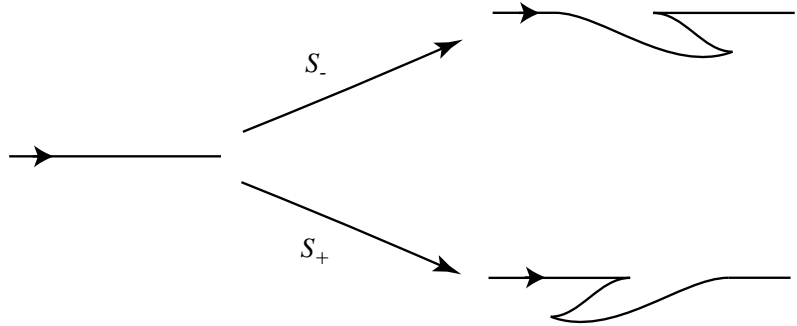}}
\caption{Positive and negative stabilizations.} 
\label{stab}	
\end{figure}

\n
{\bf Stabilization:}
Given an oriented Legendrian knot $L$, its {\em positive stabilization}
(resp.\ {\em negative stabilization}) $S_+(L)$ (resp.\
$S_-(L)$) is an operation that decreases $tb$ by adding a zigzag
in the front projection as in Figure~\ref{stab}.

\s
We have $tb(S_\pm(L)) = tb(L) - 1$ and $r(S_\pm (L))= r(L) \pm 1.$

\begin{HW}
Prove that the stabilization operation is well-defined (independent of
the location where the zigzag is added).
\end{HW}

The following theorem of Eliashberg-Fraser~\cite{EF} enumerates all the
Legendrian unknots:

\begin{thm}[Eliashberg-Fraser]\label{unknot}
Legendrian unknots in the standard contact $\R^3$ (or $S^3$) are
completely determined by $tb$ and $r$. 
\end{thm}

\begin{figure}[ht]
\epsfysize=1.7in 
\centerline{\epsfbox{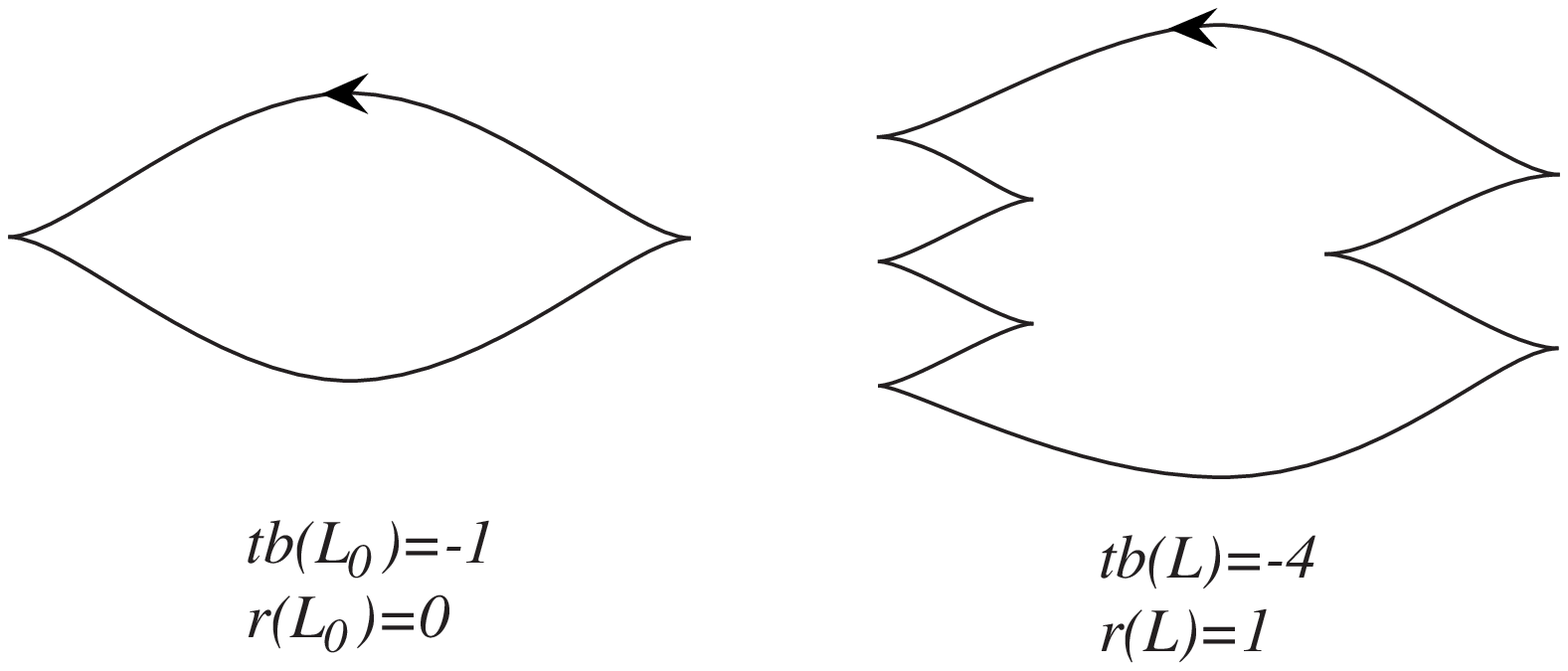}}
\caption{Legendrian unknots in the front projection.} 
\label{newunknots}	
\end{figure}

In fact, all the Legendrian unknots are stabilizations
$S^{k_1}_+S^{k_2}_-(L_0)$ of the unique maximal $tb$ Legendrian unknot
$L_0$ with $tb(L_0)=-1$ and $r(L_0)=0$, given on the left-hand side of
Figure~\ref{newunknots}.  The right-hand picture is $S^2_+ S^1_-(L_0)$. 

For an oriented Legendrian knot in $\R^3$ or $S^3$, the topological knot type,
the Thurston-Bennequin invariant, and the rotation number are called
the {\em classical} invariants.  Although Legendrian unknots are
completely determined by their classical invariants according to
Theorem~\ref{unknot}, Legendrian knots in general are not completely
classified by the classical invariants.  One way of distinguishing two
Legendrian knots with the same classical invariants is through {\em
  contact homology}. (See \cite{Ch, EGH} for more details.)

\subsection{Tight vs. overtwisted}

In the 1970's, Lutz \cite{Lu} and Martinet \cite{Ma} proved the following:

\begin{thm}[Lutz, Martinet]
Let $M$ be a closed oriented 3-manifold, $Dist(M)$ be the set of smooth
2-plane field distributions on $M$, and $Cont(M)$ be the set of smooth
contact 2-plane field distributions on $M$.  Then 
$$\pi_0(Cont(M))\rightarrow \pi_0(Dist(M))$$
is surjective.
\end{thm}

\begin{proof}[Strategy of Proof]$\mbox{}$

\be
\item Start with a 2-plane field $\xi$.  Take a fine enough
  triangulation of $M$ so that on each 3-simplex $\Delta$, $\xi$ is
  close to a linear foliation by planes. 
\item It is easy to homotop $\xi$ near the 2-skeleton so it becomes
  contact.  Now we have an extension problem to the interior of each
  3-simplex. 
\item Insert a Lutz tube.  A {\em Lutz tube} is a contact structure on
  $S^1\times D^2$ (with cylindrical coordinates $(z,r,\theta)$, where
  $D^2=\{(r,\theta)|r\leq  1\}$) given by the 1-form  
$$\alpha=\cos(2\pi r) dz+r\sin(2\pi r)d\theta.$$
\ee
\end{proof}

\begin{HW}
Think about how to use a Lutz tube (``perform a Lutz twist") to
finish the construction.  Keep in mind that the homotopy class of the
2-plane field needs to be preserved.
%Use Reeb components (see
%Section~\ref{fol}) to prove an analogous theorem for codimension 1
%foliations in 3-manifolds. 
\end{HW}

Having introduced Lutz twists, we can now write down more contact
structures on $\R^3$: 

\s\n
{\bf Example $\mbox{\bf{1}}_R$:} $(\R^3,\zeta_R)$, where $\R^3$ has
cylindrical coordinates $(r,\theta,z)$, $R$ is a positive real number,
and $\zeta_R$ is given by 
$\alpha_R=\cos f_R(r)dz + r\sin f_R(r)d\theta$.  Here $f_R(r)$ is a
function with positive derivative 
satisfying $f_R(r)=r$ near $r=0$ and $\lim_{r\rightarrow
  +\infty}f_R(r)=R$.   

\begin{HW}
Show that $(\R^3,\xi_0)\simeq (\R^3,\zeta_R)$ for all $R\leq \pi$. 
\end{HW}

However, we have the following key result of Bennequin~\cite{Be}:

\begin{thm}[Bennequin]
$(\R^3,\xi_0)\not \simeq (\R^3,\zeta_R)$ if $R>\pi$.
\end{thm}

\begin{figure}[ht]
\epsfysize=2in 
\centerline{\epsfbox{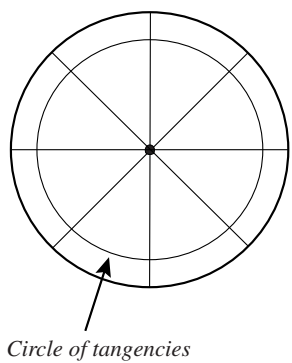}}
\caption{An overtwisted disk $D$. (Precisely speaking, the disk should
  end at the circle of tangencies.)  The straight lines represent the
  singular (characteristic) foliation that $\xi\cap TD$ traces on $D$,
  and the circle is the set of points where $\xi=TD$.  There is also
  an elliptic tangency at the center.} 	
\label{OT} 
\end{figure}

The distinguishing feature is the existence of an {\em overtwisted
  (OT) disk}, i.e., an embedded disk $D\subset (M,\xi)$ such that
  $\xi_p=T_pD$ at all $p\in \bdry D$.  A typical OT disk looks like
$\{pt\}\times D^2$ in the Lutz tube $S^1\times D^2$ described above
(also see Figure~\ref{OT}).  While it is not hard to see that
$(\R^3,\zeta_R)$ has OT disks if $R>\pi$, what Bennequin proved was
that $(\R^3,\xi_0)$ contains no OT disks.   It turns out that the
  existence of an OT disk is equivalent to the existence of a
  Legendrian unknot $L$ with $tb(L)=0$. 

\begin{HW}[Hard]
Try to prove that $(\R^3,\xi_0)$ has no overtwisted disks.
\end{HW}

It is not an exaggeration to say that modern contact geometry has its
beginnings in Bennequin's theorem.  There is a dichotomy in the world of
contact structures, those that contain OT disks (called {\em
  overtwisted} contact structures) and those that do not (called {\em
  tight} contact structures).  In view of Theorem~\ref{Pfaff}, every
contact structure is {\em locally tight}, and therefore the question
of overtwistedness is a global one.  

The following is an important inequality for knots in tight contact
manifolds. 

\begin{thm}[Bennequin inequality]
Let $L$ be nullhomologous Legendrian knot in a tight $(M,\xi)$.
If $\Sigma$ is a Seifert surface for $L$ with Euler
characteristic $\chi(\Sigma)$, then 
$$tb(L)\pm r(L)\leq -\chi(\Sigma).$$
\end{thm}

\subsection{Classification of contact structures}

When discussing the classification of contact structures, it is
important to keep in mind the following theorem: 

\begin{thm}[Gray]
Let $\xi_t$, $t\in[0,1]$, be a 1-parameter family of contact
structures on a closed manifold $M$.  Then there is a 1-parameter
family of diffeomorphisms $\varphi_t$ such that $\varphi_0=id$ and
$\varphi^*_t\xi_t=\xi_0$. 
\end{thm}

In other words, a {\em homotopy} of contact structures gives rise to a
contact {\em isotopy}.

The {\em overtwisted} classification (on closed 3-manifolds) was shown
by Eliashberg 
\cite{El3} to be essentially the same as the homotopy classification
of 2-plane fields.  (The result is quite striking, especially when
contrasted with the tight classification on $T^3$ below.)

\begin{thm}[Eliashberg]
Let $M$ be a closed oriented 3-manifold, and $Cont^{OT}(M)\subset
Dist(M)$ be the overtwisted 2-plane field distributions.  Then  
$$\pi_0(Cont^{OT}(M))\simeq \pi_0(Dist(M)).$$ 
\end{thm}

On the other hand, tight contact structures tend to reflect the
underlying topology of the manifold, and are more difficult to
understand.  The goal of this article is to introduce techniques
which enable us to better understand tight contact structures.  In the
meantime, we list a couple of examples: 

\be
\item $S^3$.  Eliashberg~\cite{El2} proved that there is a unique
  tight contact structure up to isotopy.  It is the one given in Example~3. 
\item $T^3$.  Giroux~\cite{Gi2} and Kanda~\cite{Ka} independently
  proved that (a) every tight contact structure is isomorphic to some
  $\xi_n$ and (b) $(T^3,\xi_m)\not\simeq (T^3,\xi_n)$ if $m\not=n$. 
\ee

\begin{HW}
Try to prove that $(T^3,\xi_m)\not\simeq (T^3,\xi_n)$ if $m\not=n$.
\end{HW}

In Section~\ref{three} we will give a classification of
tight contact structures for the lens spaces $L(p,q)$.

\subsection{A criterion for tightness}\label{filling}

A contact structure $(M,\xi)$ is {\em symplectically fillable} if
there exists a compact symplectic 4-manifold $(X,\omega)$ 
such that $\bdry X=M$ and $\omega|_\xi>0$.  $(X,\omega)$ is said to be
a {\em symplectic filling} of $(M,\xi)$. (Technically speaking, what
we are calling ``symplectically fillable'' is usually called
``weakly symplectically fillable'', but since we have no need of such
taxonomy in this article, we will stick to ``symplectically fillable''
or even just ``fillable''.  For more information, refer to \cite{EH}.)

\begin{HW}
Show that $(S^3,\xi)$ in Example~3 is symplectically fillable. 
\end{HW}

\begin{HW} 
Show $(T^3,\xi_n)$ in Example~2 is symplectically fillable.
(Hint: first modify $\alpha_n\mapsto dz+t\alpha_n$ with $t$ small.) 
\end{HW}

A powerful general method for producing tight contact structures is
the following theorem of Gromov and Eliashberg \cite{El1, Gr}:

\begin{thm}[Gromov-Eliashberg]
A symplectically fillable contact structure is tight.
\end{thm}

It immediately follows from the symplectic filling theorem that the
standard $(S^3,\xi)$ from Example~3 and the contact structures
$(T^3,\xi_n)$ from Example~2 are tight. 

Symplectic filling is a 4-dimensional way of checking whether
$(M,\xi)$ is tight.  We will discuss other methods (including a purely
3-dimensional one) of proving tightness in Section~\ref{four}.

\subsection{Relationship with foliation theory}\label{fol}

Foliations are the other type of locally homogeneous 2-plane field
distributions.  The following table is a brief list of analogous
objects from both worlds (note that the analogies are not precise): 

\s\s
\centerline{
\begin{tabular}{|c|c|}
\hline
Foliations & Contact Structures\\
\hline
$\alpha\wedge d\alpha=0$ & $\alpha\wedge d\alpha>0$\\
integrable & nonintegrable \\
\hline
$\alpha=dz$ & $\alpha=dz-ydx$\\
Frobenius & Pfaff\\
\hline
Reeb components & Overtwisted disks\\
\hline
Taut & Tight\\
\hline
\end{tabular}}

\begin{figure}[ht]
\epsfysize=2in 
\centerline{\epsfbox{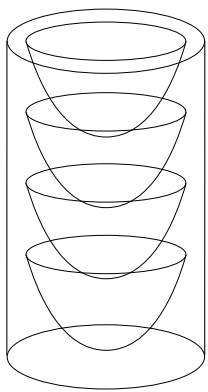}}
\caption{A Reeb component.  Here the top and bottom are identified.} 	
\label{Reeb} 
\end{figure}

\s\s\n
A (rank 2) {\em foliation} $\xi$ is an integrable 2-plane field
distribution, i.e., locally given as the kernel of a 1-form $\alpha$ with
$\alpha\wedge d\alpha=0$.  According to Frobenius' theorem, $\xi$ can
locally be written as the kernel of $\alpha=dz$.  The world of
foliations also breaks up into the 
topologically significant {\em taut} foliations (i.e., foliations for
which there is a closed transversal curve through each leaf), and the
foliations with {\em generalized Reeb components}, which exist on
every 3-manifold.  A generalized Reeb component is a compact submanifold
$N\subset M$ whose boundary $\bdry N$ is a union of torus leaves, and
such that there are no transversal arcs which begin and end on $\bdry
N$.  The primary example of a generalized Reeb component
is a {\em Reeb component}, i.e., a foliation of the
solid torus $S^1\times D^2$ whose boundary $S^1\times S^1$ is a
leaf and whose interior is foliated by planes as in Figure~\ref{Reeb}.

The following is a key theorem which allows us to transfer
information from foliation theory to contact geometry.

\begin{thm} [Eliashberg-Thurston]
Let $M$ be a closed, oriented 3-manifold $\not=S^1\times S^2$.  Then
every taut foliation admits a $C^0$-small perturbation into a tight contact
structure. 
\end{thm}

For a thorough treatment of the relationship with foliation theory,
see~\cite{ET}.  In Section~\ref{decomp}, we will discuss one aspect,
namely the relationship with Gabai's {\em sutured manifold} theory.

\section{Convex surfaces}

In this section, we investigate embedded surfaces $\Sigma$ in the
contact manifold $(M,\xi)$.  The principal notion is that of {\em
  convexity}.  For the time being, $\xi$ may be tight or overtwisted.

\subsection{Characteristic foliations}\label{charfol}

Before discussing convexity, we first examine how $\xi$
traces a singular line field on an embedded surface $\Sigma$.

\begin{defn}
The {\em characteristic foliation} $\Sigma_\xi$ is the singular
foliation induced on $\Sigma$ from $\xi$, where
$\Sigma_\xi(p)=\xi_p\cap T_p\Sigma$.  The {\em singular points} (or
tangencies) are points $p\in \Sigma$ where $\xi_p=T_p\Sigma$. 
\end{defn}

\begin{lemma}\label{generic}
A $C^\infty$-generic characteristic foliation $\Sigma_\xi$ is of
{\em Morse-Smale type}, i.e., satisfies the following:
\be
\item the singularities and closed orbits are {\em dynamically
  hyperbolic}, i.e, hyperbolic in the dynamical systems sense,
\item there are no saddle-saddle connections, and
\item every point $p\in \Sigma$ limits to some isolated singularity or 
  closed orbit in forward time and likewise in backward time.
\ee
\end{lemma}

The proof of Lemma~\ref{generic} uses the fact that a $C^\infty$-small
perturbation of $\xi$ is still contact.  We choose the perturbation of $\xi$ to be
compactly supported near $\Sigma$, and hence the isotopy in Gray's theorem is
compactly supported near $\Sigma$.  Therefore, generic properties of
1-forms (in particular the Morse-Smale condition) are satisfied. 

\begin{HW}
Show that if $\alpha$ is a contact 1-form and $\beta$ is any 1-form,
then $\alpha+t\beta$ is contact for sufficiently small $t$. 
\end{HW}

There are two types of dynamically hyperbolic singularities: 
{\em elliptic} and {\em hyperbolic} (not in the dynamical systems
sense).  Choose coordinates $(x,y)$ on $\Sigma$ and let the origin be 
the singular point.  If we write $\alpha= dz+fdx +gdy$, then
$X=g{\bdry\over \bdry x}-f{\bdry\over \bdry y}$ is a vector field for the
characteristic foliation near the origin.  If the determinant of the
matrix
$$\begin{pmatrix}
{\bdry g\over \bdry x} & {\bdry g\over \bdry y}\\
-{\bdry f\over \bdry x} & -{\bdry f\over \bdry y}
\end{pmatrix}$$
is positive (resp.\ negative), then the singular point is elliptic
(resp.\ hyperbolic).  An example of an elliptic singularity is
$\alpha=dz+(xdy-ydx)$, and an example of a hyperbolic singularity is
$\alpha=dz+(2xdy+ydx)$. 

Next we discuss signs.  Assume $\Sigma$ and $\xi$ are both oriented.
Then a singular point $p$ is {\em positive} (resp.\ {\em negative}) if
$T_p\Sigma$ and $\xi_p$ have the same orientation (resp.\ opposite
orientations). 

\begin{claim}
The characteristic foliation $\Sigma_\xi$ is oriented.  
\end{claim}

We use the convention that positive elliptic points are sources and
negative elliptic points are sinks.  If $p$ is a nonsingular point of
a leaf $L$ of the characteristic foliation, then we choose $v\in T_pL$
so that $(v,n)$ is an oriented basis for $T_p\Sigma$.  Here $n\in
T_p\Sigma$ is an oriented normal vector to $\xi_p$. 

\s\n
{\bf Examples of characteristic foliations:}
\s
\be
\item Consider $S^2=\{x^2+y^2+z^2=1\}\subset (\R^3,\zeta_{\pi/2})$.
  Then $S^2$ will have two singular points, the positive elliptic
  point $(0,0,1)$ and the negative elliptic point $(0,0,-1)$, and the
  leaves spiral downward from $(0,0,1)$ to $(0,0,-1)$. 

\item An example of an overtwisted disk $D$ is one which has a
  positive elliptic point at the center and radial leaves emanating
  from the center, such that $\bdry D$ is a circle of singularities.
  Often in the literature one sees overtwisted disks whose boundary is
  transverse to $\xi$ and whose leaves emanating from the center
  spiral towards the limit cycle $\bdry D$. (Strictly speaking,
  such a $D$ with a limit cycle is not an OT disk according to our
  definition, but can easily be modified to fit our definition.)
\ee

\s
The importance of the characteristic foliation $\Sigma_\xi$ comes from
the following proposition: 

\begin{prop}
Let $\xi_0$ and $\xi_1$ be two contact structures which induce the
same characteristic foliation on $\Sigma$.  Then there is an isotopy
$\varphi_t$, $t\in[0,1]$, rel $\Sigma$, with $\varphi_0=id$ and
$(\varphi_1)_*\xi_0=\xi_1$. 
\end{prop}

\subsection{Convexity}

The notion of a {\em convex surface}, introduced by Giroux in \cite{Gi1} and
extended to the case of a compact surface with Legendrian boundary by
Kanda in \cite{Ka}, is the key ingredient in the cut-and-paste theory of contact
structures. 

\begin{defn}
A properly embedded oriented surface $\Sigma$ is {\em convex} if there
exists a contact vector field $v\pitchfork \Sigma$.  Here, a {\em
  contact vector field} is a vector field whose corresponding flow
preserves the contact structure $\xi$.  In this article we assume that
our convex surfaces are either {\em closed} or {\em compact with
  Legendrian boundary}.   
\end{defn}

If $\Sigma=\Sigma\times\{0\}$ is convex, then there is an invariant
neighborhood $\Sigma\times [-\varepsilon,\varepsilon]\subset M$.  We
usually assume that $v$ agrees with the normal orientation to $\Sigma$. 

\s\n
{\bf Properties of convex surfaces:}

\be
\item A $C^\infty$-generic closed embedded surface $\Sigma$ is convex.
  This is because an embedded surface $\Sigma$ with a Morse-Smale
  characteristic foliation is convex. (The
  same is almost true for compact surfaces with Legendrian boundary,
  but more care is needed along the boundary.) 

\item To a convex surface $\Sigma$ we may associate a multicurve
  (i.e., a properly embedded (smooth) 1-manifold, possibly
  disconnected and possibly with boundary) 
$$\Gamma_\Sigma=\{x\in \Sigma | v(x)\in \xi_x\},$$
called the {\em dividing set}.  It satisfies the following:

\be
\item $\Gamma_\Sigma\pitchfork \Sigma_\xi$.

\item The isotopy class of $\Gamma_\Sigma$ does not depend on the
  choice of $v$. 

\item $\Sigma\setminus \Gamma_\Sigma = R_+(\Gamma_\Sigma)\sqcup
  R_-(\Gamma_\Sigma)$, where $R_+(\Gamma_\Sigma)\subset \Sigma$
  (resp.\ $R_-(\Gamma_\Sigma)$) is the set of points 
  $x$ where the normal orientation to $\Sigma$ given by $v(x)$ agrees
  with (resp.\ is opposite to) the normal orientation to $\xi_x$. 
\ee

\ee

\begin{rmk}
We may think of $\Gamma_\Sigma$ as the set of points where $\xi\perp
\Sigma$, where $\perp$ is measured with respect to $v$. 
\end{rmk}

Write $\#\Gamma_\Sigma$ for the number of connected components of
$\Gamma_\Sigma$. 

\s
\begin{figure}[ht]
\epsfysize=1.5in 
\centerline{\epsfbox{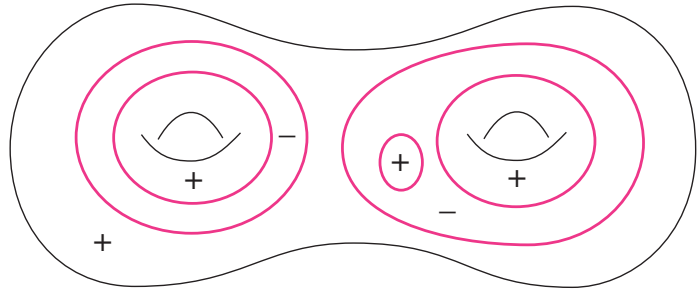}}
\caption{A sample dividing set.} 	
\label{dividing} 
\end{figure}

\s
The usefulness of the dividing set $\Gamma_\Sigma$ comes from the following:

\begin{thm}[Giroux's Flexibility Theorem]
Assume $\Sigma$ is convex with characteristic foliation $\Sigma_\xi$,
contact vector field $v$, and dividing set $\Gamma_\Sigma$.  Let
$\mathcal{F}$ be another singular foliation on $\Sigma$ which is {\em
  adapted to} $\Gamma_\Sigma$ (i.e., there is a contact structure
$\xi'$ in a neighborhood of $\Sigma$ such that
$\Sigma_{\xi'}=\mathcal{F}$ and $\Gamma_\Sigma$ is also a dividing set
for $\xi'$).  Then there is an isotopy $\varphi_t$, $t\in [0,1]$, of
$\Sigma$ in $(M,\xi)$ such that: 
\be
\item $\varphi_0=id$ and $\varphi_t|_{\Gamma_\Sigma}=id$ for all $t$.
\item $\varphi_t(\Sigma)\pitchfork v$ for all $t$.
\item $\varphi_1(\Sigma)$ has characteristic foliation $\mathcal{F}$.
\ee
\end{thm}

In essence, $\Gamma_\Sigma$ encodes ALL of the essential
contact-topological information in a neighborhood of $\Sigma$.
Therefore, having discussed characteristic foliations in
Section~\ref{charfol}, we may proceed to discard them and simply
remember the dividing set.

\begin{HW}
Prove Giroux Flexibility.
\end{HW}

\n
{\bf Examples on $T^2$:}  There are two common characteristic
foliations on $T^2$.
\be
\item {\em Nonsingular Morse-Smale}.  This is when the characteristic
  foliation is nonsingular and has exactly $2n$ closed orbits, $n$ of
  which are sources (repelling periodic orbits) and the other $n$ are
  sinks (attracting periodic orbits).  $\Gamma_{T^2}$ 
  consists of $2n$ closed curves parallel to the closed orbits.  Each
  dividing curve lies inbetween two periodic orbits.  

\item {\em Standard form}.  An example is $x=const$ inside
  $(T^3,\xi_n)$.  The torus is fibered by closed Legendrian fibers,
  called {\em ruling curves}, and the singular set consists of $2n$
  closed curves, called {\em Legendrian divides}.  The $2n$ curves of
  $\Gamma_{T^2}$ lie between the Legendrian divides.  
\ee

\begin{HW}
Find an explicit example of a $T^2$ inside a contact manifold with
nonsingular Morse-Smale characteristic foliation. 
\end{HW}

\begin{figure}[ht]
\epsfysize=2in 
\centerline{\epsfbox{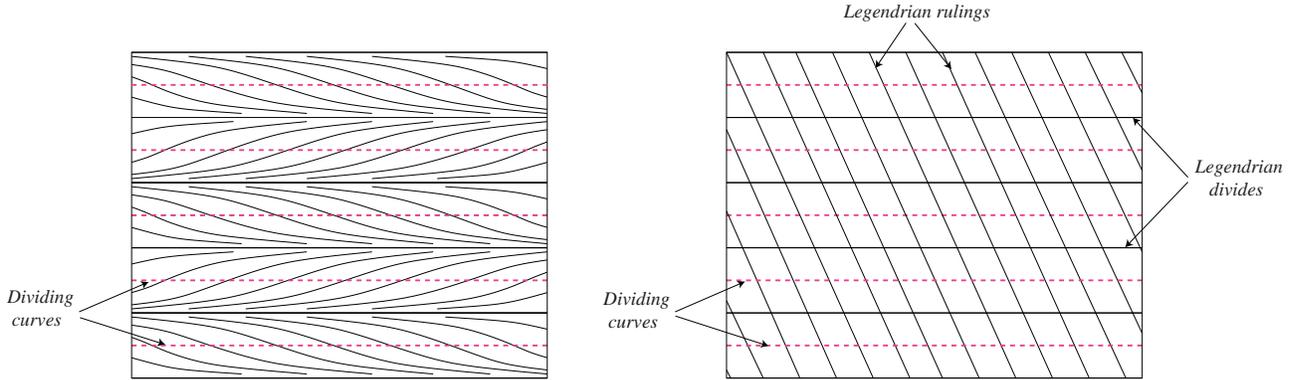}}
\caption{The left-hand side is a torus with nonsingular Morse-Smale
  characteristic foliation.  The right-hand side is a torus in
  standard form.  Here the sides are identified and the top and bottom
  are identified.} 	 
\end{figure}

What Giroux Flexibility tells us is that it is easy to switch between
the two types of characteristic foliations -- nonsingular Morse-Smale
and standard form. The following corollary of Giroux Flexibility is a
crucial ingredient in the cut-and-paste theory of contact structures.

\begin{cor}[Legendrian Realization Principle, abbreviated LeRP]
Let $\Sigma$ be a convex surface and $C$ be a multicurve on $\Sigma$.
Assume $C\pitchfork \Gamma_\Sigma$ and 
$C$ is {\em nonisolating}, i.e., each connected component of
$\Sigma\setminus C$ nontrivially intersects $\Gamma_\Sigma$. 
Then there is an isotopy (as in the Giroux Flexibility Theorem) such
that $\varphi_1(C)$ is Legendrian. 
\end{cor}

\begin{HW}
Try to prove LeRP, assuming Giroux Flexibility.
\end{HW}

\begin{rmk}
$C$ may have extraneous intersections with $\Gamma_\Sigma$, i.e., the
  actual number of intersections $\#(C\cap \Gamma_\Sigma)$ is allowed
  to be larger than the geometric intersection number. 
\end{rmk}

\n
{\bf Fact:} If $C$ is a Legendrian curve on the convex surface
$\Sigma$, then the twisting number $t(C,\Sigma)$ relative to the
framing from $\Sigma$ is $-{1\over 2}\#(C\cap
\Gamma_\Sigma)$. Here $\#(\cdot)$ represents cardinality, not geometric
intersection. 

\s
Now we present the criterion (see \cite{Gi3}) for determining when a
convex surface has a tight neighborhood. 

\begin{prop}[Giroux's Criterion]
A convex surface $\Sigma\not=S^2$ has a tight neighborhood if and only
if $\Gamma_\Sigma$ has no homotopically trivial dividing curves.  If
$\Sigma=S^2$, then there is a tight neighborhood if and only if
$\#\Gamma_\Sigma=1$.  
\end{prop}

\begin{HW}\label{19}
Prove that if $\Gamma_\Sigma$ has a homotopically trivial dividing
curve, then there exists an overtwisted disk in a neighborhood of
$\Sigma$, provided we are not in the situation where $\Sigma=S^2$ and
$\#\Gamma_\Sigma=1$.  (Hint: use LeRP, together with a trick when 
$\Gamma_\Sigma$ has no other components besides the homotopically
trivial curve.) 
\end{HW}

The ``only if'' direction in Giroux's Criterion follows from HW~\ref{19}.  The
``if'' direction follows from constructing an explicit model inside a tight
3-ball or gluing (for the latter, see \cite{Co1}).

Suppose that $(M,\xi)$ is tight.  If $\Sigma=S^2$
is a convex surface in $(M,\xi)$, then $\Gamma_\Sigma$ is unique up to
isotopy, consisting on one (homotopically trivial) circle.  If
$\Sigma=T^2$ is convex, then it consists of $2n$ parallel,
homotopically essential curves.  Therefore $\Gamma_{T^2}$ is
determined by $\#\Gamma_{T^2}$ and the slope, once a trivialization
$T^2\simeq \R^2/\Z^2$ is fixed.

\subsection{Convex decomposition theory}\label{decomp}

The reader may have already noticed certain similarities between
convex surfaces and the theory of sutured manifolds due to
Gabai~\cite{Ga}.  

\begin{defn}
A {\em sutured manifold} $(M,\Gamma)$ consists of the following data:
\be
\item $M$ is a compact, oriented, irreducible 3-manifold; each
  component of $M$ has nonempty boundary, 
\item $\Gamma$ is a multicurve on $\bdry M$ which has nonempty
  intersection with each component of $\bdry M$, and
\item $\Gamma$ divides $\bdry M$ into positive and negative regions,
  whose sign changes every time $\Gamma$ is crossed.  We write
  $\bdry M\setminus \Gamma=R_+(\Gamma) \sqcup R_-(\Gamma)$. 
\ee
Here, a 3-manifold $M$ is {\em irreducible} if every embedded 2-sphere
$S^2$ bounds a 3-ball $B^3$. 
\end{defn} 

Note that our definition of a sutured manifold, chosen to simplify the
exposition in this paper, is slightly different from that of
Gabai~\cite{Ga}. 

\begin{defn} 
Let $S$ be a compact oriented surface with connected components
$S_1,\dots,S_n$.  The {\em Thurston norm} of $S$ is:
$$ x(S)= \sum_{\mbox{$i$ such that $\chi(S_i)<0$}} |\chi(S_i)|.$$
\end{defn}

\begin{defn}
A sutured manifold $(M,\Gamma)$ is {\em taut} if $R_\pm(\Gamma)$  are
incompressible in $M$ and minimize the Thurston norm in 
$H_2(M,\Gamma)$.  Here, a surface $S\subset M$ is {\em incompressible}
if for every embedded disk $D\subset M$ with $D\cap S=\bdry D$, there
is a disk $D'\subset S$ such that $\bdry D=\bdry D'$.
\end{defn}

Roughly speaking, $(M,\Gamma)$ is taut if $R_\pm(\Gamma)$ attain the
minimum genus amongst all the embedded representatives in the relative
homology class $H_2(M,\Gamma)$.

We have the following theorem which gives the equivalence between
tightness and tautness in the case of a manifold with boundary (see
\cite{HKM1}): 

\begin{thm}[Kazez-Mati\'c-Honda]\label{maintheorem} 
Let $(M,\Gamma)$ be a sutured manifold.  Then the following are equivalent: 
\be
\item $(M,\Gamma)$ is taut.
\item $(M,\Gamma)$ carries a taut foliation.
\item $(M,\Gamma)$ carries a universally tight contact structure.
\item $(M,\Gamma)$ carries a tight contact structure.
\ee
\end{thm}

A contact structure $\xi$ on $M$ is {\em carried by $(M,\Gamma)$} if
$\bdry M$ is a convex surface for $\xi$ with dividing set $\Gamma$.  A
transversely oriented foliation $\xi$ on $M$ is {\em carried
  by $(M,\Gamma)$} if there exists a 
thickening of $\Gamma$ to a union $\gamma\subset \bdry M$ of annuli,
so that $\bdry M\setminus \gamma$ is a union of leaves of $\xi$, $\xi$
is transverse to $\gamma$, and the orientations of $R_\pm(\Gamma)$ and
$\xi$ agree.  (Strictly speaking, in this case $M$ is a manifold with
corners.)   A tight contact structure is {\em universally tight} if
it remains tight when pulled back to the universal cover of $M$. 

In the rest of this section, we explain how {\em sutured manifold
  decompositions} have an analog in the contact world, namely the
theory of {\em convex decompositions}.  Using it we outline the proof
of (1)$\Rightarrow$(4).

\begin{figure}[ht]
\epsfysize=1.2in 
\centerline{\epsfbox{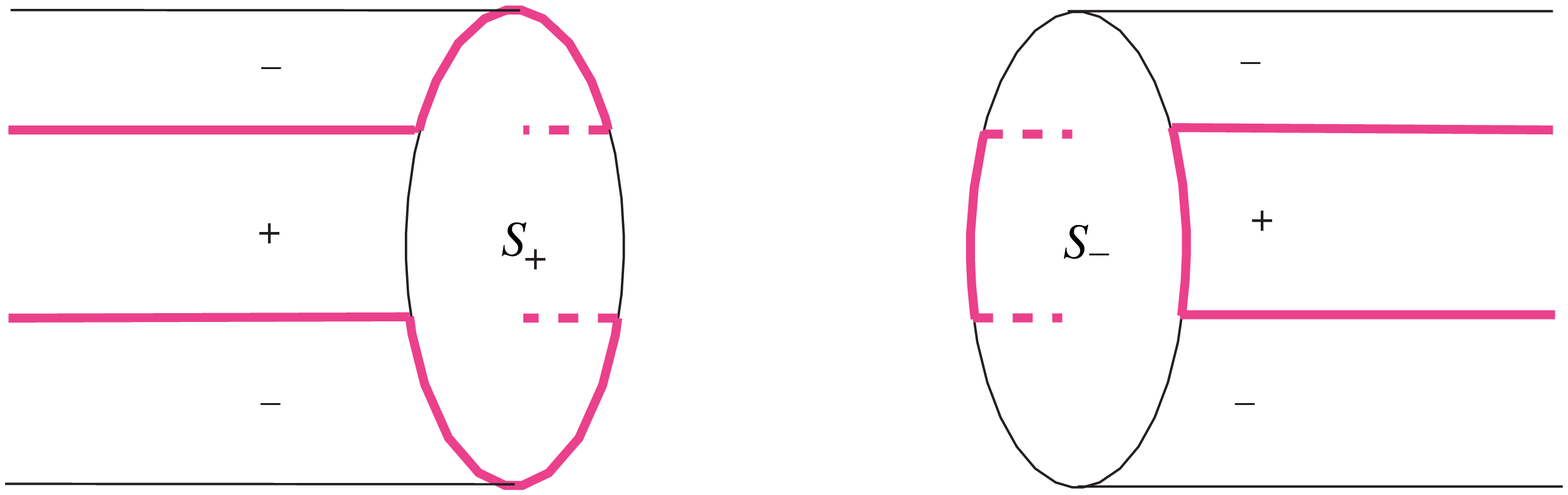}}
\caption{} 	 
\label{suturedsplitting2} 
\end{figure}

\begin{defn}
Let $S$ be an oriented, properly embedded surface in $(M,\Gamma)$
which intersects $\Gamma$ transversely.  Then a {\em sutured manifold
  splitting}
$(M,\Gamma)\stackrel{S}\sa (M',\Gamma')$ is given
as follows (see Figure~\ref{suturedsplitting2} for an illustration): 
Define $M'=M\setminus S$, and let $S_+$ (resp.\ $S_-$) be the copy of
$S$ on $\bdry 
M'$ where the orientation inherited from $S$ and the outward normal
agree (are opposite).  Then set $R_\pm
(\Gamma')=(R_\pm(\Gamma)\setminus $S$)\cup S_\pm$.  The new suture
$\Gamma'$ forms the boundary between the regions $R_+(\Gamma')$ and
$R_-(\Gamma')$.
\end{defn}

\s

\s
A sutured manifold $(M,\Gamma)$ is {\em decomposable}, if there is a
sequence of sutured manifold splittings:
$$(M,\Gamma)\stackrel{S_1}\sa (M_1,\Gamma_1)\stackrel{S_2}\sa 
\dots\stackrel{S_n}\sa (M_n,\Gamma_n)=\sqcup (B^3,S^1).$$

\s\n
Gabai, in \cite{Ga}, proved the following theorem:

\begin{thm}[Gabai] $\mbox{}$
\be
\item (Decomposition) If $(M,\Gamma)$ is taut, then it is decomposable.
\item (Reconstruction) Given a sutured manifold decomposition, we can
  backtrack and construct a taut foliation which is carried by $(M,\Gamma)$.
\ee
\end{thm}

Now, in the contact category, we choose a dividing set $\Gamma_S$ so
that every component of $\Gamma_S$ is {\em $\bdry$-parallel},  i.e.,
cuts off a half-disk of $S$ which does not intersect any other
component of $\Gamma_S$.  Such a dividing set $\Gamma_S$ is also
called {\em $\bdry$-parallel}.  

\begin{figure}[ht]
\epsfysize=1in 
\centerline{\epsfbox{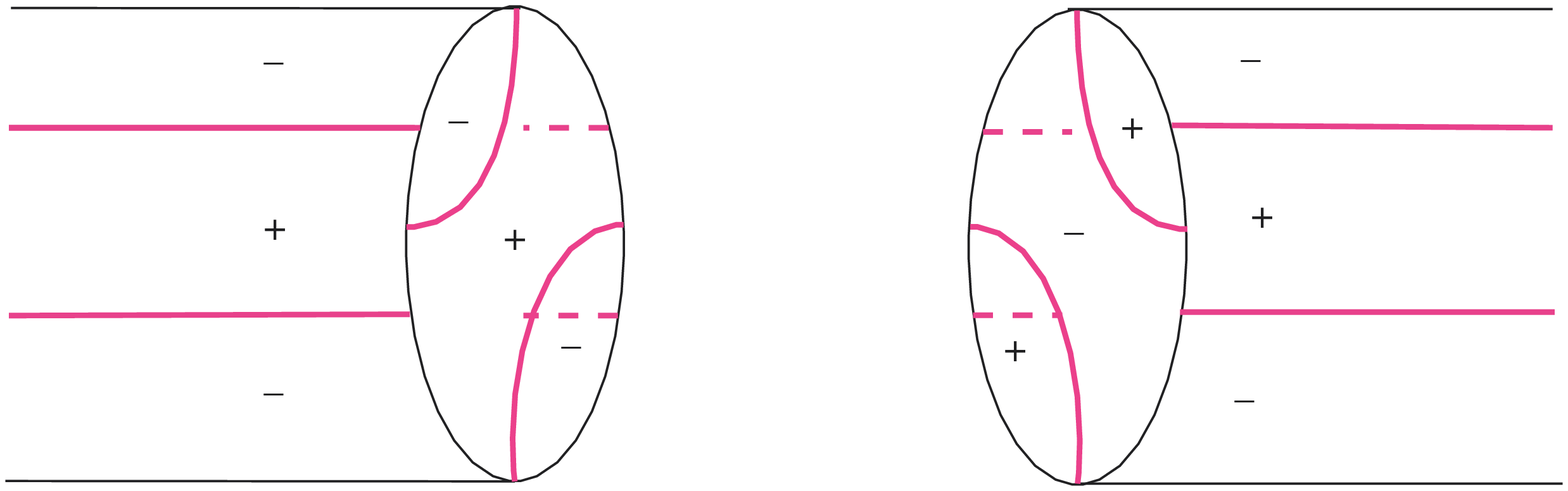}}
\caption{} 	
\label{rounding1} 
\end{figure}

\begin{figure}[ht]
\epsfysize=1in 
\centerline{\epsfbox{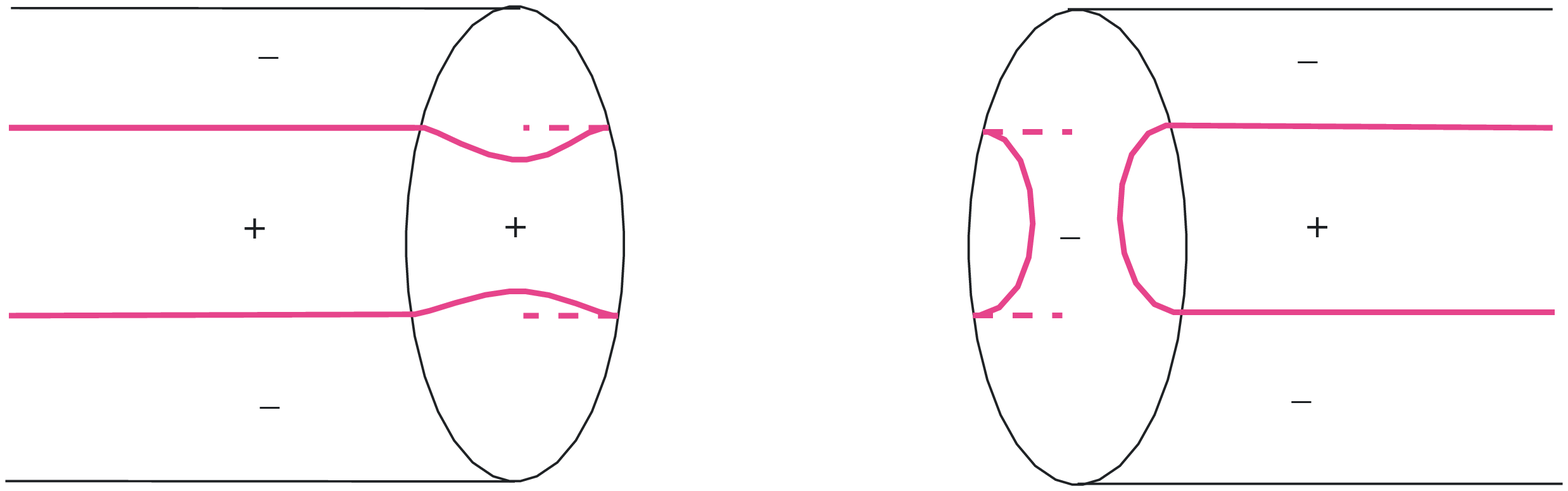}}
\caption{} 	
\label{rounding2} 
\end{figure}

If there is an invariant contact structure defined in a neighborhood
of $\bdry M$ with dividing set $\Gamma=\Gamma_{\bdry M}$, then by an
application of LeRP, we may take $\bdry S$ to be Legendrian.  (There
are some exceptional cases, but we will not worry about them here.)
Extend the contact structure to be an invariant contact structure in a
neighborhood of $S$ with $\bdry$-parallel dividing set $\Gamma_S$.
Now, if we cut $M$ along $S$, we obtain a manifold with corners.   To
smooth the corners, we apply {\em edge-rounding}.  This is given in 
Figures~\ref{rounding1} and \ref{rounding2}.  Figure~\ref{rounding1}
gives the surface $S$ before rounding, and Figure~\ref{rounding2}
after rounding.  Notice that we may think of $S$ as a lid of a
jar, and the edge-rounding operation as twisting to close the jar.  
\s

\s
\begin{HW}
Explain why edge-rounding works as in Figures~\ref{rounding1} and
\ref{rounding2}. 
\end{HW}

Observe that the dividing set in Figure~\ref{rounding2} is isotopic to
the sutures in Figure~\ref{suturedsplitting2}.
Therefore, given a sutured manifold splitting
$(M,\Gamma)\stackrel{S}\sa (M',\Gamma'),$
there is a corresponding convex splitting
$(M,\Gamma)\stackrel{(S,\Gamma_S)}\sa (M',\Gamma'),$ with a
$\bdry$-parallel dividing set $\Gamma_S$. 
Using the decomposition theorem of Gabai, if $(M,\Gamma)$ is taut,
then there exists a convex decomposition:
$$(M,\Gamma)\stackrel{(S_1,\Gamma_{S_1})}\sa (M_1,\Gamma_1)
\stackrel{(S_2,\Gamma_{S_2})}\sa  
\dots\stackrel{(S_n,\Gamma_{S_n})}\sa 
(M_n,\Gamma_n)=\sqcup (B^3,S^1).$$

We now work backwards, starting with the following theorem of
Eliashberg \cite{El2}:

\begin{thm}[Eliashberg]
Fix a characteristic foliation $\mathcal{F}$ adapted to $\Gamma_{\bdry
  B^3}=S^1$.  Then there is a unique tight contact structure on $B^3$
up to isotopy relative to $\bdry B^3$. 
\end{thm}

The following gluing theorem of Colin \cite{Co1} allows us to
inductively build a universally tight contact structure carried by
$(M,\Gamma)$. 

\begin{thm}[Colin] \label{u-tight}
Let $\Sigma$ be an incompressible surface with
$\bdry\Sigma\not=\emptyset$.  If $\Gamma_\Sigma$ is $\bdry$-parallel
and $(M\setminus \Sigma,\xi|_{M\setminus \Sigma})$ is universally
tight, then $(M,\xi)$ is also universally tight.   
\end{thm}

This theorem and other similar theorems will be discussed in
Section~\ref{four}. 

Theorem~\ref{maintheorem} is a refinement, in the case of manifolds
with boundary, of the following theorem:

\begin{thm}[Gabai-Eliashberg-Thurston]
Let $M$ be a oriented, closed, irreducible 3-manifold with
$H_2(M;\Z)\not=0$.  Then $M$ carries a universally tight contact
structure. 
\end{thm}

The Gabai-Eliashberg-Thurston theorem was originally proved in two
parts: Gabai \cite{Ga} proved that such an $M$ carries a taut
foliation, and Eliashberg-Thurston \cite{ET} proved that the taut
foliation can be perturbed into a universally tight contact
structure.  There is also an alternate, purely 3-dimensional method
for proving this theorem \cite{HKM2,HKM3,HKM4}.

\section{Bypasses}

In this section, we introduce the other chief ingredient in the
cut-and-paste theory of tight contact structures: the {\em bypass}.
As a surface is isotoped inside the ambient tight contact manifold
$(M,\xi)$, the dividing set changes in discrete units, and the
fundamental unit of change is effected by the bypass.  Bypasses would
be quite useless if they were difficult to find.  For the cases we
examine in Section~\ref{three}, namely solid tori, $T^2\times I$, and
lens spaces, they can be found relatively easily by examining the next
step in the Haken hierarchy.  This will be explained in
Section~\ref{find}.  For more information on bypasses, refer to \cite{H1}.

\subsection{Definition and examples}

\begin{defn}
Let $\Sigma$ be a convex surface and $\alpha$ be a Legendrian arc in
$\Sigma$ which intersects $\Gamma_\Sigma$ in three points
$p_1,p_2,p_3$, where $p_1$ and $p_3$ are endpoints of $\alpha$.  A
{\em bypass half-disk} is a convex half-disk $D$ with Legendrian
boundary, where $D\cap \Sigma=\alpha$ and $tb(\bdry D) =-1$.  $\alpha$
is called the {\em arc of attachment} of the bypass, and $D$ is said
to be a bypass {\em along} $\alpha$ or $\Sigma$.
\end{defn}

\begin{figure}[ht]
\epsfysize=2in 
\centerline{\epsfbox{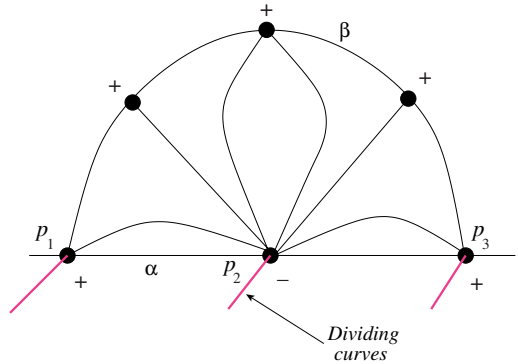}}
\caption{A bypass.} 	
%\label{} 
\end{figure}

\begin{rmk}
Most bypasses do not come for free.  Finding a bypass is equivalent to
raising the twisting number (or Thurston-Bennequin invariant) by 1.
Although it is easy to lower the twisting number by attaching
``zigzags'' in a front projection, raising the twisting number is
usually a nontrivial operation.  
\end{rmk}

\begin{figure}[ht]
\epsfysize=1.8in 
\centerline{\epsfbox{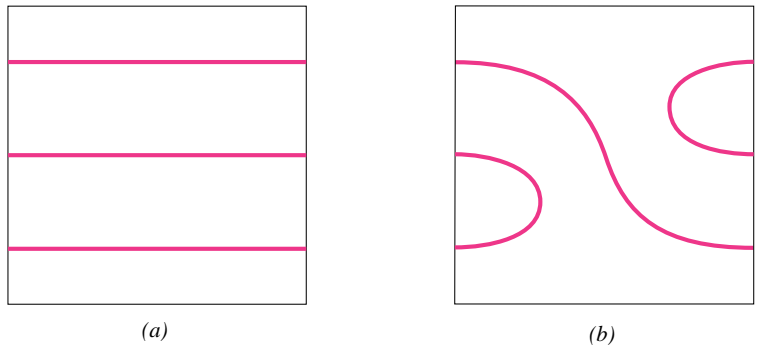}}
\caption{The effect of attaching a bypass from the
  front. $\Gamma_\Sigma$ is (a) and $\Gamma_{\Sigma'}$ is (b).}  
\label{bypassmove} 
\end{figure}

\begin{lemma}[Bypass Attachment Lemma]
Let $D$ be a bypass for $\Sigma$.  If $\Sigma$ is isotoped across $D$,
then we obtain a new convex surface $\Sigma'$ whose dividing set is
obtained from $\Gamma_\Sigma$ via the move in
Figure~\ref{bypassmove}. 
\end{lemma}

Note that this is reasonable because a bypass attachment increases the
twisting number along the arc of attachment by 1.

\begin{figure}[ht]
\epsfysize=3.5in 
\centerline{\epsfbox{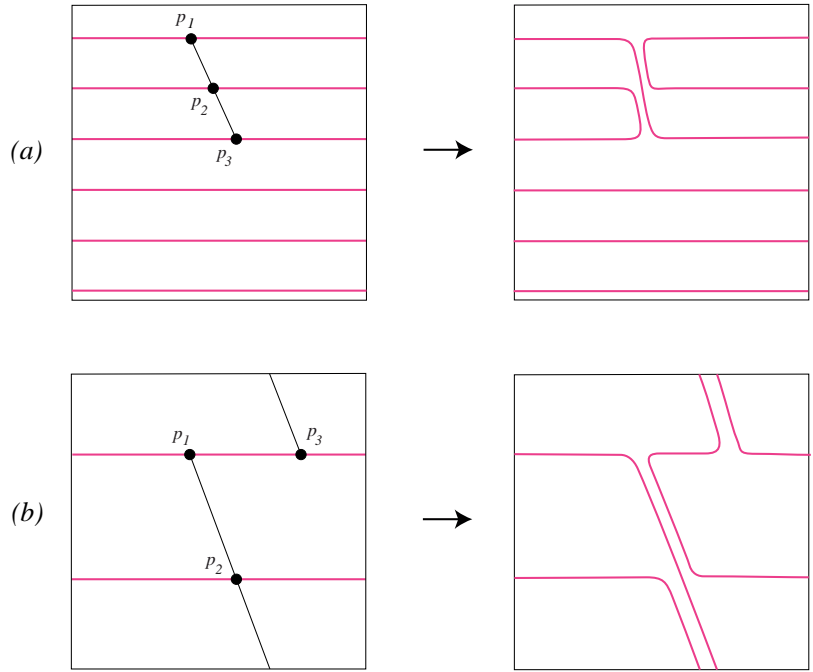}}
\caption{Possible bypasses on tori.} 	
\label{torus} 
\end{figure}

\begin{figure}[ht]
\epsfysize=1.6in 
\centerline{\epsfbox{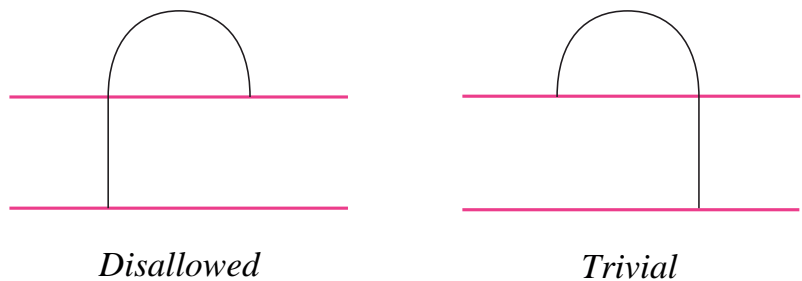}}
\caption{A disallowed bypass attachment and a trivial bypass
  attachment.}  	
\label{disallowed} 
\end{figure}

\s\n
{\bf Example:} $T^2$.  Let us enumerate the possible bypass
attachments -- see Figure~\ref{torus}.  (a) is the case where
$\#\Gamma_{T^2}=2n>2$, and the bypass reduced $\#\Gamma$ by two, while
keeping the slope fixed.  (b) is the case where $\#\Gamma_{T^2}=2$,
and the slope is modified.  In addition, there also are trivial and
disallowed moves, which are moves locally given in
Figure~\ref{disallowed}. It turns out that the trivial move always
exists inside a tight contact manifold, whereas the disallowed move
can never exist inside a tight contact manifold. 

\begin{HW}
Is there a bypass attachment which increases $\#\Gamma$?
\end{HW}

\begin{figure}[ht]
\epsfysize=2in 
\centerline{\epsfbox{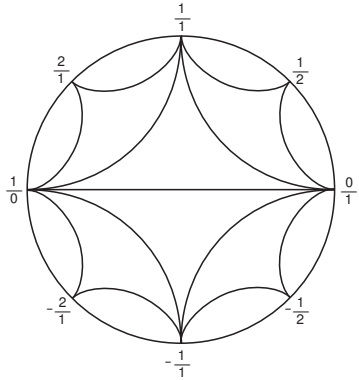}}
\caption{The Farey tessellation. The spacing between vertices are not 
  drawn to scale.} 	
\label{farey} 
\end{figure}

\n
{\bf Intrinsic interpretation:}
Observe that, in case (b), the bypass move is equivalent to performing
a {\em positive Dehn twist} along a particular curve.  We can
therefore reformulate this bypass move and give an {\em intrinsic
  interpretation} in terms of the Farey tessellation of the
hyperbolic unit disk $\H$ (Figure~\ref{farey}).  The set of vertices
of the Farey tessellation is $\Q\cup\{\infty\}$ on 
$\bdry \H$.  (More precisely, fix a fractional linear transformation $f$
from the upper half-plane model of hyperbolic space to the unit disk
model $\H$.  Then the set of vertices is the image of
$\Q\cup\{\infty\}$ under $f$.)  There is a unique edge between
${p\over q}$ and ${p'\over q'}$ if and only if the
corresponding shortest integer vectors form an integral basis for
$\Z^2$. (The edge is usually taken to be a geodesic in $\H$.)

\begin{prop}\label{counterclockwise}
Let $s=\mbox{slope}(\Gamma_{T^2})$.  If a bypass is attached along a
closed Legendrian curve of slope $s'$, then the resulting slope $s''$
is obtained as follows: Let $(s',s)\subset \bdry \H$ be the
counterclockwise interval from $s'$ to $s$.  Then $s''$ is the point
on $(s',s)$ which is closest to $s'$ and has an edge to $s$.
\end{prop}

\begin{figure}[ht]
\epsfysize=2in 
\centerline{\epsfbox{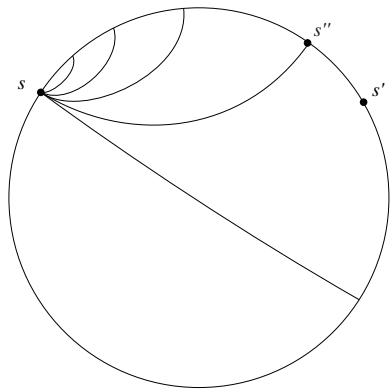}}
\caption{Intrinsic interpretation of the bypass attachment.}
\label{farey-bypass}
\end{figure}

See Figure~\ref{farey-bypass} for an illustration.

\begin{HW}
Prove Proposition~\ref{counterclockwise}.
\end{HW}

\subsection{Finding bypasses}\label{find}
We now explain how to find bypasses.  Let $M$ be a closed manifold
and $\Sigma\subset M$ be a closed surface.  In order to find a bypass
along $\Sigma$, we consider $M\setminus \Sigma$.  Let $S\subset M\setminus \Sigma$ 
be an incompressible surface with nonempty boundary, for example the
next cutting surface in the Haken hierarchy.  Under mild conditions on
$\bdry S$, we can take $S$ to be a convex surface with nonempty
Legendrian boundary.

\begin{lemma}\label{bp}
Suppose that $\Gamma_S$ has a $\bdry$-parallel component and either
$S\not=D^2$ or else if $S=D^2$ then $tb(\bdry S)<-1$.   Then there
exists a bypass along $\bdry S$ and hence along $\Sigma$. 
\end{lemma}

\begin{proof}
Draw an arc $\delta'\subset S$ so that $\delta'$ cuts off a half-disk
with only the $\bdry$-parallel arc $\delta$ on it.  The condition on
$S$ is needed to ensure that we can use LeRP to find a Legendrian arc
$\delta''$.  The half-disk cut off by $\delta''$ (and containing a
copy of $\delta$) is the bypass for $\Sigma$. 
\end{proof}

\begin{cor} \label{disk}
Let $S=D^2$ be a convex disk with Legendrian boundary so that
$tb(\bdry S)<-1$.  Then there exists a bypass along $\bdry S$.
\end{cor}

Corollary~\ref{disk} follows from Lemma~\ref{bp}, by observing that
all components of $\Gamma_{D^2}$ cut off half-disks of $D^2$ and that
a $\bdry$-parallel component is simply an outermost arc of
$\Gamma_{D^2}$. 

\begin{rmk}
Corollary~\ref{disk} does not work when $tb(\bdry D)=-1$.
\end{rmk}
 
Similarly, we can prove the following:

\begin{cor}[Imbalance Principle]
Let $S=S^1\times[0,1]$ be a convex annulus.  If
$t(S^1\times\{1\},\mathcal{F}_S)<t(S^1\times\{0\},\mathcal{F}_S)$,  
then there is a $\bdry$-parallel arc and hence a bypass along
$S^1\times\{1\}$.  Here $\mathcal{F}_S$ is the framing induced from
the surface $S$.
\end{cor}

\begin{figure}[ht]
\epsfysize=2in 
\centerline{\epsfbox{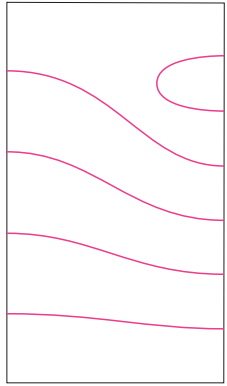}}
\caption{One possible dividing set for the annulus.  Here the top and
  the bottom are identified.}  	
\label{annulus} 
\end{figure}

Figure~\ref{annulus} gives an example of a convex annulus with
$t(S^1\times\{1\},\mathcal{F}_S)<t(S^1\times\{0\},\mathcal{F}_S)$.
There is necessarily a bypass along $S^1\times\{1\}$.

\section{Classification of tight contact structures on lens
  spaces}\label{three} 

As an illustration of the technology introduced in the previous two
sections, we give a complete classification of tight contact
structures on the lens spaces $L(p,q)$.  This classification was
obtained independently by Giroux~\cite{Gi2} and Honda~\cite{H1};
partial results had been obtained previously by Etnyre~\cite{Et2}.  In
this article, we follow the method of \cite{H1}.

\subsection{The standard neighborhood of a Legendrian curve}

Consider a (closed) Legendrian curve $L$ with
$t(L,\mathcal{F})=-n<0$, $n\in \Z^+$.  (Pick some framing
$\mathcal{F}$ for which the twisting number is negative.)  Then a {\em
  standard neighborhood} $S^1\times D^2=\R/\Z \times \{x^2+y^2\leq
\varepsilon\}$ (with coordinates $z,x,y$) of the Legendrian curve
$L=S^1\times \{ (0,0)\}$ is given by  
$$\alpha = \sin (2\pi nz)dx +\cos(2\pi n z)dy,$$
and satisfies the following:

\be
\item $T^2= \bdry (S^1\times D^2)$ is convex.
\item $\#\Gamma_{T^2}=2$.
\item $\mbox{slope}(\Gamma_{T^2})=-{1\over n}$, if the meridian has
  zero slope and the longitude given by $x=y=const$ has slope
  $\infty$. 
\ee

The following is due to Kanda~\cite{Ka} and Makar-Limanov~\cite{ML1}.

\begin{prop}[Kanda, Makar-Limanov]
Given a solid torus $S^1\times D^2$ and boundary conditions (1), (2),
(3), there exists a unique tight contact structure on $S^1\times D^2$
up to isotopy rel boundary, provided we have fixed a characteristic
foliation $\mathcal{F}$ adapted to $\Gamma_{\bdry(S^1\times T^2)}$. 
\end{prop}

\begin{rmk}
The precise characteristic foliation is irrelevant in view of Giroux
Flexibility. 
\end{rmk}

\begin{proof} $\mbox{}$

\be
\item Let $L\subset T^2$ be a curve which bounds the meridian
  $D$. Using LeRP, realize it as a Legendrian curve with $tb(L)=-1$.   

\item Using the genericity of convex surfaces, realize the surface $D$
  with $\bdry D=L$ as a convex surface with Legendrian boundary.
  Since $tb(L)=-1$, there is only one possibility for
  $\Gamma_{D}$, up to isotopy.   

\item Next, using Giroux Flexibility, fix some characteristic
  foliation on $D$ adapted to $\Gamma_{D}$.  Note that any two tight
  contact structures on $S^1\times D$ with boundary condition
  $\mathcal{F}$ can be isotoped to agree on $T^2\cup D$. 

\item The rest is a 3-ball $B^3$.  Use Eliashberg's uniqueness theorem
  for tight contact structures on $B^3$.
\ee 
\end{proof}

\begin{HW}
Try to prove Eliashberg's theorem, using convex surfaces.
\end{HW}

\subsection{Lens spaces}

Let $p>q>0$ be relatively prime integers.  The {\em lens space}
$L(p,q)$ is obtained by gluing $V_1=S^1\times D^2$ and $V_2=S^1\times
D^2$ together via $A: \bdry V_2\stackrel{\sim}{\rightarrow} \bdry
V_1$, where $A=\begin{pmatrix} -q & q' \\ p & -p' 
\end{pmatrix} \in -SL(2,\Z)$.  Here we are making an oriented
identification $\bdry V_i\simeq \R^2/\Z^2$, where the meridian of
$V_i$ is mapped to $\pm (1,0)$, and some chosen longitude is mapped to
$\pm (0,1)$. 

\s\n
{\bf Continued fractions:}  Let $-{p\over q}$ have a continued
fraction expansion 
$$-{p\over q} = r_0 -{1\over r_1-{1\over r_2\dots -{1\over r_k}}},$$
where $r_i\leq -2$.   

\s\n
{\bf Example:} $-{14\over 5} = -3-{1\over -5}$.  We write $-{14\over
  5} \leftrightarrow (-3,-5)$.

\begin{thm}[Giroux, Honda]\label{lens}
On $L(p,q)$, there are exactly $|(r_0+1)(r_1+1)\dots(r_k+1)|$ tight
contact structures up to isotopy.  They are all {\em holomorphically
  fillable}. 
\end{thm}

A surgery presentation for $L(p,q)$ is given as follows:

\begin{figure}[ht]
\epsfysize=1.3in 
\centerline{\epsfbox{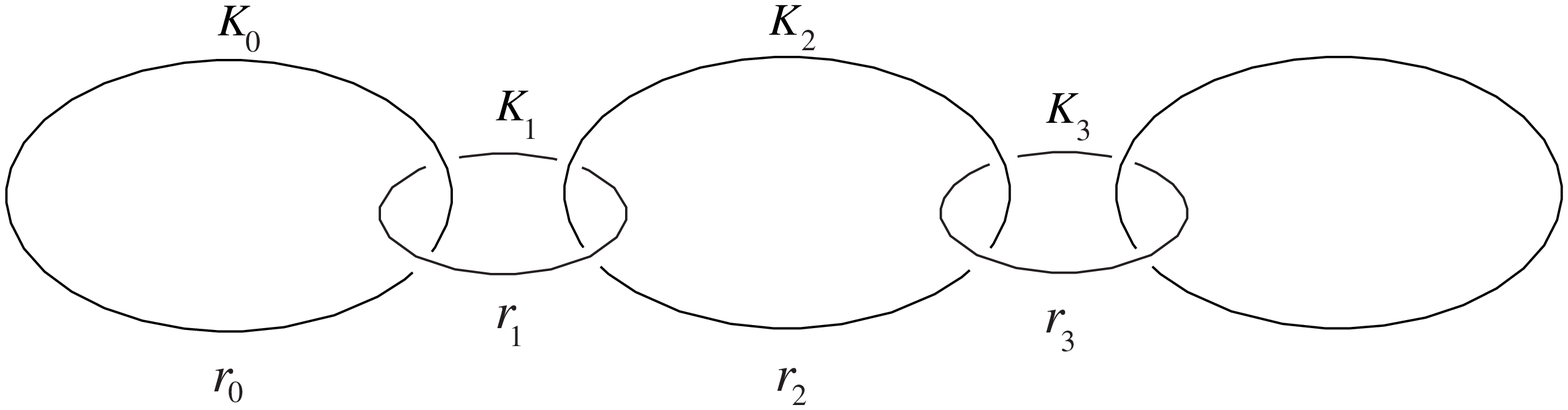}}
\caption{} 	
\label{link} 
\end{figure}

\s\n
{\bf Legendrian surgery:}  Given a Legendrian knot $K=K_0$ or link
$L=\sqcup_{i=0}^k K_i$ in a contact manifold $(M,\xi)$, we can
perform a surgery along the $K_i$ with coefficient $tb(K_i)-1$.  At
the 4-dimensional level, if $M=S^3$, then we start with a Stein domain
$B^4$ with $\bdry B^4=S^3$, and attach 2-handles in a way
which makes the resulting 4-manifold $X^4$ a Stein domain (and in
particular symplectic).  The resulting contact 3-manifold $(M',\xi')$
with $\bdry X=M'$ is said to be {\em holomorphically
  fillable}. Similarly, if $(M,\xi)$ is symplectically fillable, then
$(M',\xi')$ obtained by Legendrian surgery is also symplectically
fillable.  The Stein construction was done by Eliashberg in \cite{El4}
and the symplectic construction by Weinstein~\cite{We}.

\s
Suppose $K_i$ is a Legendrian unknot with $tb(K_i)=r_i+1$
and $r(K_i)= \mbox{one of } r_i+2, r_i+4,\dots, -(r_i+2).$  There
are precisely $|r_i+1|$ choices for the rotation number
$r(K_i)$. (In fact, these are all the Legendrian unknots with
$tb(K_i)=r_i+1$ by Theorem~\ref{unknot}.) 

\begin{HW}
Show that the $|r_0+1| |r_1+1| \dots |r_k+1|$ holomorphically fillable
contact structures are distinct. 
\end{HW}

Therefore, we have the lower bound:
\begin{equation}
\# \mbox{Tight}(L(p,q))\geq |(r_0+1)(r_1+1)\dots (r_k+1)|.
\end{equation}
\n
Here $\mbox{Tight}(M)$ refers to the set of isotopy classes of tight
contact structures on $M$.
In order to prove Theorem~\ref{lens}, it remains to show the reverse
inequality.

\subsection{Solid tori}  
We now consider tight contact structures on
the solid torus $S^1\times D^2$ with the following conditions on the
boundary $T=S^1\times D^2$: 
\be
\item $\#\Gamma_T=2$.
\item $\mbox{slope}(\Gamma_T)=-{p\over q}$, where $-\infty<-{p\over
  q}\leq -1$.  (After performing Dehn twists, we can normalize the
  slope as such.) 
\item The fixed characteristic foliation $\mathcal{F}$ is adapted to
  $\Gamma_T$. 
\ee

\begin{thm}\label{solid torus}
There are exactly $|(r_0+1)(r_1+1)\dots (r_{k-1}+1) r_k|$ tight
contact structures on $S^1\times D^2$ with this boundary condition. 
\end{thm}

\n
{\bf Step 1:} In this step we factor $S^1\times D^2$ into a union of
$T^2\times I$ layers and a standard neighborhood of a Legendrian
curve isotopic to the core curve of $S^1\times D^2$.   Assume
$-{p\over q}<-1$, since $-{p\over q}=-1$ has already 
been treated.   

Let $D$ be a meridional disk with $\bdry D$ Legendrian and $tb(\bdry
D)= -p<-1$.  Then by Lemma~\ref{disk} there is at least one bypass
along $\bdry D$.  Attach the bypass to $T$ from the interior and apply
the Bypass Attachment Lemma.  We obtain a convex torus $T'$ isotopic
to $T$, such that $T$ and $T'$ cobound a $T^2\times I$.  Denote
$\mbox{slope}(\Gamma_{T'})=-{p'\over q'}$. 

\begin{HW}
If $-{p\over q}\leftrightarrow (r_0,r_1,\dots,r_{k-1},r_k)$, then
$-{p'\over q'}\leftrightarrow (r_0,r_1,\dots,r_{k-1},r_k+1)$. 
\end{HW}

We successively peel off $T^2\times I$ layers according to the Farey
tessellation.  The sequence of slopes is given by the continued
fraction expansion, or, equivalently, by the shortest sequence of
counterclockwise arcs in the Farey tessellation from $-{p\over q}$ to
$-1$.  Once slope $-1$ is reached, $S^1\times D^2$ with boundary slope $-1$
is the standard neighborhood of a Legendrian core curve with twisting
number $-1$ (with respect to the fibration induced from the $S^1$-fibers
$S^1\times\{pt\}$). 

\s\n
{\bf Step 2:} (Analysis of each $T^2\times I$ layer)

\s\n
{\bf Fact:}  Consider $T^2\times [0,1]$ with convex boundary
conditions $\#\Gamma_0=\#\Gamma_1=2$,  $s_0=\infty$, and $s_1=0$.
Here we write $\Gamma_i=\Gamma_{T^2\times\{i\}}$ and
$s_i=\mbox{slope}(\Gamma_i)$.  (More invariantly, the shortest
integers corresponding $s_0,s_1$ form an integral basis for $\Z^2$.)
Then there are exactly two tight contact structures (up to isotopy
rel boundary) which are {\em minimally twisting}, i.e., every convex
torus $T'$ isotopic to $T^2\times\{i\}$ has
$\mbox{slope}(\Gamma_{T'})$ in the interval $(0,+\infty)$.  They are
distinguished by the Poincar\'e duals of the {\em relative half-Euler
  class}, which are computed to be  $\pm ((1,0)-(0,1))\in
H_1(T^2\times[0,1];\Z)$. 
We call these $T^2\times[0,1]$ layers {\em basic slices}. 

\s
The proof of the fact will be omitted, but one of the key elements in
the proof is the following lemma: 

\begin{HW}
Prove, using the Imbalance Principle, that for any tight contact
structure on $T^2\times[0,1]$ with boundary slopes $s_0\not=s_1$ and
any rational slope $s$ in the interval $(s_1,s_0)$, there exists a
convex surface $T'\subset T^2\times[0,1]$, which is parallel to
$T^2\times\{pt\}$ and has slope $s$.  Here, if
$s_0< s_1$, $(s_1,s_0)$ means $(s_1,+\infty] \cup [-\infty,s_0)$. 
\end{HW}

\s\n
{\bf Step 3:}  (Shuffling) Consider the example of the solid torus
where $-{p\over q}= -{14\over 5}$.  We have the following
factorization: 
$$\begin{array}{ccc}
-{14\over 5} & \leftrightarrow & (-3,-5)\\
-{11\over 4} & \leftrightarrow & (-3,-4)\\
-{8\over 3} & \leftrightarrow & (-3,-3)\\
-{5\over 2} &\leftrightarrow & (-3,-2)\\
- 2 & \leftrightarrow & (-3,-1)=(-2)\\
-1 &\leftrightarrow & (-1)
\end{array}$$
We group the basic slices into {\em continued fraction blocks}.  Each
block consists of all the slopes whose continued fraction
representations are of the same length.  In the example, we have two
blocks: slope $-{14\over 5}$ to $-2$, and slope $-2$ to $-1$.  All the
relative half-Euler classes of the basic slices in the first block are $\pm
(-1,3)$; for the second block, they are $\pm (0, 1)$.  Therefore, a
naive upper bound for the number of tight contact structures would be
$2$ to the power $\#(\mbox{basic slices})$.  

A closer inspection however reveals that we may {\em shuffle} basic
slices which are in the same continued fraction block.  More
precisely, if $T^2\times[0,2]$ admits a factoring into basic slices
$T^2\times[0,1]$ and $T^2\times[1,2]$ with relative half-Euler classes
$(a,b)$ and $-(a,b)$, then it also admits a factoring into basic
slices where the relative half-Euler classes are $-(a,b)$ and $(a,b)$,
i.e., the order is reversed.  

Shuffling is (more or less) equivalent to the following proposition:

\begin{lemma}\label{commute}
Let $L$ be a Legendrian knot.  Then $S_+S_-(L)=S_-S_+(L)$.
\end{lemma}

\begin{HW}
Prove Lemma~\ref{commute}.  (Observe that the ambient contact manifold
is irrelevant and that the commutation can be done in a standard
tubular neighborhood of $L$.)
\end{HW}

Returning to the example at hand, the first continued fraction block
has at most $|-5| = 4+1$ tight contact structures (distinguished by the relative
half-Euler class), and the second has at most $|-3+1| = 2$ tight
contact structures.  We compute $\#\mbox{Tight} \leq 2\cdot 5$.    

In general, for the solid torus with slope $-{p\over q}\leftrightarrow
(r_0,r_1,\dots, r_k)$ we have: 

\begin{equation} \label{upperbound}
\#\mbox{Tight} \leq |(r_0+1)(r_1+1)\dots (r_{k-1}+1) r_k|.
\end{equation}

\subsection{Completion of the proof of Theorems~\ref{lens} and
  \ref{solid torus}}  

We prove the following, which instantaneously completes the proof of
both theorems.  
\begin{equation}\label{lensupper}
\# \mbox{Tight}(L(p,q))\leq |(r_0+1)(r_1+1)\dots (r_k+1)|.
\end{equation}
Recall that on $\bdry V_1$, the meridian of $V_2$ has slope $-{p\over
  q}\leftrightarrow (r_0,r_1,\dots, r_{k-1},r_k)$.  First, take a
Legendrian curve $\gamma$ isotopic to the core curve of $V_2$ with
largest twisting number.  (Such a Legendrian curve exists, since any
  closed curve admits a $C^0$-small approximation by a Legendrian
  curve; the upper bound exists by the Thurston-Bennequin inequality.)
  We may assume $V_2$ is the standard 
neighborhood of $\gamma$; the tight contact structure on $V_2$ is
then unique up to isotopy.   Next, $\mbox{slope}(\Gamma_{\bdry
  V_1})=-{p'\over q'}\leftrightarrow (r_0,\dots, r_{k-1}, r_k+1)$, and
we have already computed the upper bound for $\#\mbox{Tight} (V_2)$ to
be $|(r_0+1)\dots(r_{k-1}+1)(r_k+1)|$ by Equation~\ref{upperbound}.
  This completes the proof of Equation~\ref{lensupper} and hence of
  Theorems~\ref{lens} and \ref{solid torus}.

\begin{oq}
Give a complete classification of tight contact structures on $T^2\times
[0,1]$ when $\#\Gamma_{T^2\times\{i\}}>2$, $i=0,1$.  (Contrary to what
is claimed in \cite{H1}, the general answer is not yet known.)
\end{oq}

\section{Gluing}\label{four}

There are three general methods for proving tightness:
\be
\item symplectic filling,
\item gauge theory (in particular Heegaard Floer homology), and
\item gluing (state traversal).
\ee
Symplectic filling was already discussed in Section~\ref{filling}.
We briefly explain the relationship between contact structures and the
{\em Heegaard Floer homology} of Ozsvath and Szabo \cite{OSz1,OSz2}.
To an oriented closed 3-manifold $M$ one can assign a {\em Heegaard
  Floer homology group} $\widehat{HF}(M)$, constructed out of the
Heegaard decomposition of $M$.  In \cite{OSz3}, Ozsvath and Szabo
assigned a class $c(\xi)\in \widehat{HF}(-M)$ to every 
contact structure $(M,\xi)$ (tight or overtwisted).  This was done via
the work of Giroux~\cite{Gi4} in which it was shown that every contact
structure $(M,\xi)$ corresponds to an equivalence class of open book
decompositions of $M$ (and hence an equivalence class of fibered
knots).  Lisca and Stipsicz \cite{LS3} showed that large families of
contact structures are tight (but not fillable) by showing that their Heegaard
Floer homology class is nonzero.  The Heegaard Floer homology approach appears
to be very promising at the time of the writing of this article.

In this section we focus on the last technique, namely gluing.  Many
of the key ideas in gluing were introduced by Colin~\cite{Co1,Co2} and
Makar-Limanov~\cite{ML2}, and subsequently enhanced by Honda \cite{H2}
who combined them with the bypass technology.

Let us start by asking the following question:

\begin{q}
Let $\Sigma$ be a convex surface in $(M,\xi)$.  If $(M\setminus
\Sigma, \xi|_{M\setminus \Sigma})$ is tight, then is $(M,\xi)$ tight? 
\end{q}

\n
{\bf Answer:} This is usually not true.  Our goal is to understand to
what extent it is true. 

\begin{HW}
Give an example of an overtwisted $T^2\times [0,1]$ which is tight
when restricted to $T^2\times[0,{1\over2 }]$ and to $T^2\times[{1\over
    2},1]$. 
\end{HW}

\subsection{Basic examples with trivial state transitions}$\mbox{ }$

\s\n
{\bf Example A:} (Colin \cite{Co2}, Makar-Limanov \cite{ML2})  Suppose
$\Sigma=S^2$.  If $(M\setminus \Sigma,\xi|_{M\setminus \Sigma})$ is
tight, then $(M,\xi)$ is tight. 

\begin{proof}
Recall that there is only one possibility for $\Gamma_{S^2}$ inside a
tight contact manifold.  We argue by contradiction.  Suppose there is
an OT disk $D\subset M$.  A priori, the OT disk $D$ can intersect
$\Sigma$ in a very complicated manner.  We obtain a contradiction as
follows: 

\be
\item Isotop $\Sigma$ to $\Sigma'$ so that $\Sigma'\cap D=\emptyset$.
\item {\em Discretize} the isotopy 
$$\Sigma_0=\Sigma \rightarrow \Sigma_1\rightarrow \dots \rightarrow
  \Sigma_n=\Sigma',$$ 
so that each step is obtained by attaching a {\em bypass}.
\item If $(M\setminus \Sigma_i,\xi|_{M\setminus \Sigma_i})$ is tight,
  then $\Gamma_{\Sigma_i} = \Gamma_{\Sigma_{i+1}}= S^1$ and the bypass
  must be {\em trivial}.  Hence, 
$$(M\setminus \Sigma_i,\xi|_{M\setminus \Sigma_i})\simeq (M\setminus
  \Sigma_{i+1},\xi|_{M\setminus \Sigma_{i+1}}).$$ 
\ee
We have proved inductively that $(M\setminus \Sigma',\xi|_{(M\setminus
  \Sigma')})$ is tight, a contradiction. 
\end{proof}

More generally, one can prove:

\begin{thm}[Colin~\cite{Co2}]
If $M=M_1\# M_2$, then
$$ \mbox{Tight}(M)\simeq \mbox{Tight}(M_1)\times \mbox{Tight}(M_2).$$
\end{thm}

\begin{HW}
Classify tight contact structures on $S^1\times S^2$.
\end{HW}

\n
{\bf Example B:} (Colin~\cite{Co1}) If $\Sigma=D^2$ and
$\Gamma_\Sigma$ is {\em $\bdry$-parallel}, then $(M\setminus
\Sigma,\xi|_{M\setminus \Sigma})$ tight $\Rightarrow (M,\xi)$ tight.    
 
\s\n
{\bf Example C:} (Colin~\cite{Co1}) Let $\Sigma$ be an incompressible
surface with $\bdry\Sigma\not=\emptyset$.  If $\Gamma_\Sigma$ is
$\bdry$-parallel and $(M\setminus \Sigma,\xi|_{M\setminus \Sigma})$ is
universally tight, then $(M,\xi)$ is universally tight.   (This is
Theorem~\ref{u-tight} above.)

\begin{q}
In Example C, does $(M\setminus \Sigma,\xi|_{M\setminus \Sigma})$
tight imply $(M,\xi)$ tight? In other words, can universal tightness
be avoided? 
\end{q}

All of the above examples can be characterized by the fact that the
{\em state transitions} are trivial.  However, to create more
interesting examples, we need to ``traverse all states''.

\subsection{More complicated example} $\mbox{}$

\s\n
{\bf Example D:} (Honda \cite{H2}) Let $H$ be a handlebody of genus $g$ and
$D_1,\dots,D_g$ be compressing disks so that $H\setminus (D_1\sqcup
\dots\sqcup D_g)=B^3$. Fix $\Gamma_{\bdry H}$ (and a compatible
characteristic foliation).  Note that we need $tb(D_i)\leq -1$, since
otherwise we can find an OT disk using LeRP. 

Let $\mathcal{C}$ be the {\em configuration space}, i.e., the set of
all possible $C=(\Gamma_{D_1},\dots,\Gamma_{D_g})$, where each
$\Gamma_{D_i}$ has no closed curves.  The cardinality of $\mathcal{C}$
is {\em finite}.  If we cut $H$ along $\Sigma=D_1\cup\dots\cup
D_g$, then we obtain a 3-ball with corners. 
Given a configuration $C$, we can round the corners, as previously
explained in Section~\ref{decomp}.  Now, if $\Gamma_{\bdry
  (H\setminus\Sigma)}=S^1$ after rounding, 
then $C$ is said to be {\em potentially allowable}. 

\s\n
{\bf State transitions:}  The smallest unit of isotopy (in the contact
world) is a bypass attachment.  Therefore we examine the effect of one
bypass attachment onto $D_i$.  First we need to ascertain whether a
candidate bypass exists.   

\s\n
{\bf Criterion for existence of state transition:}
The candidate bypass exists if and only if attaching the bypass from
the interior of $B^3=H\setminus \Sigma$ does not increase
$\#\Gamma_{\bdry B^3}$. 

\s
We construct a graph $\Gamma$ with $\mathcal{C}$ as the
vertices. We assign an edge from   
$(\Gamma_{D_1},\dots,\Gamma_{D_i},\dots,\Gamma_{D_g})$ to
$(\Gamma_{D_1},\dots,\Gamma_{D_i'},\dots,\Gamma_{D_g})$ if there is a
state transition $D_i\rightarrow D_i'$ given by a single bypass
move. Note that the bypass may be from either side of $D_i$.
Then we have: 

\begin{thm}
$\mbox{Tight}(H,\Gamma_{\bdry H})$ is in 1-1 correspondence with the
  connected components of $\Gamma$, all of whose vertices $C$ are
  potentially allowable. 
\end{thm} 

\begin{HW}
Explain why $\mbox{Tight}(H,\Gamma_{\bdry H})$ is finite.
\end{HW}

\begin{rmk}
Since $\mathcal{C}$ is a finite graph, in theory we can compute
$\mbox{Tight}(H,\Gamma_{\bdry H})$ for any handlebody $H$ with a fixed
boundary $\Gamma_{\bdry H}$.  Tanya Cofer, a (former) graduate student at the
University of Georgia, has programmed this for $g=1$, and the
experiment agrees with the theoretical number from Theorem~\ref{lens}, in case
$\#\Gamma_{\bdry H}=2$ and the slope is $-{p\over q}$ with $p\leq
10$. 
\end{rmk}

\begin{HW}
Using the state transition technique, analyze tight contact structures
on $S^1\times D^2$, where $\Gamma_{T^2}$, $T^2=\bdry (S^1\times D^2)$,
satisfies the following: 
\be
\item $\#\Gamma_{T^2}=2$ and $\mbox{slope}(\Gamma_{T^2})=-2$.
\item $\#\Gamma_{T^2}=2$ and $\mbox{slope}(\Gamma_{T^2})=-3$.
\item $\#\Gamma_{T^2}=4$ and $\mbox{slope}(\Gamma_{T^2})=\infty$.
\ee
Here the slope of the meridian is $0$ and the slope of some preferred
longitude is $\infty$. 
\end{HW}

\subsection{Tightness and fillability}

We present two examples which show that the world of tight contact
structures is larger than the world of symplectically fillable contact
structures. 

\s\n
{\bf Example E:} (Honda \cite{H2}) We present a tight handlebody $H$
of genus $4$ which 
becomes OT after a Legendrian surgery.  Since Legendrian surgery
preserves fillability, the tight handlebody cannot be embedded inside
any closed fillable contact 3-manifold. 

\begin{figure}[ht]
\epsfysize=2.3in 
\centerline{\epsfbox{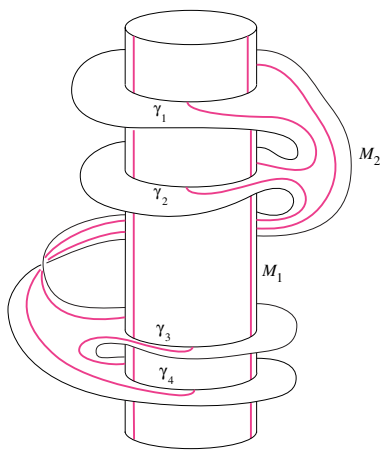}}
\caption{The top and bottom are identified.} 	
\label{handlebody} 
\end{figure}

We take the union $H=M_1\cup M_2$, where $M_1=S^1\times D^2$ is the
standard tubular neighborhood of a Legendrian curve and $M_2$ is an
$I$-invariant neighborhood of a convex disk $S$ with 4 holes.  Here
$\bdry S=\gamma-\cup_{i=1}^4 \gamma_i$ and $\Gamma_S$ consists of 4
arcs, one each from $\gamma_i$ to $\gamma_{i+1}$ ($i$ mod 4).  The
gluing is presented in Figure~\ref{handlebody}, where
$T^2=\bdry(S^1\times D^2)$ is drawn so that $\Gamma_{T^2}$ has slope
$\infty$, the $\gamma_i$ have slope $0$, and the meridian of $M_1$ has
slope $1$. 

A Legendrian surgery along the core curve of $M_1$ yields a new
meridional slope of $0$ along $T^2$, and hence allows $S$ to be
completed to an OT disk.  Using the state transition method, one can
prove that the contact structure is tight. 

\begin{HW}
Verify the tightness.
\end{HW}

\s\n
{\bf Example F:} (Etnyre-Honda~\cite{EH}) Consider the torus bundle
$M= (T^2\times[0,1])/\sim$, where $(x,1)\sim (Ax,0)$, $T^2=\R^2/\Z^2$,
and $A=\begin{pmatrix} 0 & 1 \\ -1 & 0 \end{pmatrix}.$   Let
$T^2\times[0,1]$ be a basic slice with boundary slopes $s_0=\infty$
and $s_1=0$.  The glued-up contact structure $\xi$ is proved to be 
tight using state traversal.  However, $\xi$ is not symplectically
fillable by the following contradiction argument:  

\be
\item $M$ is a Seifert fibered space over $S^2$ with Seifert
  invariants $(-{1\over 2}, {1\over 4}, {1\over 4})$.   
\item There exists a Legendrian surgery taking $(M,\xi)$ to
  $(M',\xi')$, where $M'$ is a Seifert fibered space over $S^2$ with
  invariants $(-{1\over 2}, {1\over 3}, {1\over 4})$.  Since
  Legendrian surgery preserves fillability, if $\xi$ is fillable, then
  $\xi'$ is also fillable.  
\item A theorem of Lisca~\cite{Li}, proved using Seiberg-Witten
  theory, states that there are no fillable contact structures on
  $M'$.  
\ee

\begin{rmk}
Example F was the first example of a tight contact structure which is
not fillable.  Since then, numerous other examples have been
discovered by Lisca and Stipsicz~\cite{LS1,LS2,LS3}.   
\end{rmk}

\begin{oq}
Elucidate the difference between the world of tight contact structures
and the world of fillable contact structures.
\end{oq}

\end{document}